\documentclass[reqno]{amsart}
\newcommand\presentday{{November 9, 2026}}
\usepackage{amssymb,latexsym}
\usepackage{amsmath}
\usepackage{xcolor}
\usepackage{graphicx}
\usepackage{hyperref}

\numberwithin{equation}{section}
\newtheorem{theorem}{Theorem}[section]
\newtheorem{proposition}[theorem]{Proposition}

\newtheorem{corollary}[theorem]{Corollary}

\theoremstyle{definition}
\newtheorem{remark}[theorem]{Remark}
\newtheorem{definition}[theorem]{Definition}

\newcommand{\hrefdoi}[1]{\href{https://doi.org/#1}{DOI:~#1}}
\newcommand\figwidthcoeff{0.95}

\newcommand \cleq {\mathrel{\leq_{\mathrm{co}}}}
\newcommand \eleq {\mathrel{\leq_{\mathrm{et}}}}
\DeclareMathOperator{\Coal}{Coal}

\DeclareMathOperator{\fil}{fil}
\DeclareMathOperator{\spf}{spf}
\DeclareMathOperator{\halfspf}{spf/2}
\DeclareMathOperator{\idl}{idl}
\DeclareMathOperator{\Min}{Min}
\DeclareMathOperator{\Max}{Max}
\DeclareMathOperator{\UpCov}{UpCov}
\DeclareMathOperator{\ET}{ET}

\DeclareMathOperator{\coord}{coord}
\DeclareMathOperator{\rec}{rec}
\DeclareMathOperator{\dep}{dep}
\DeclareMathOperator{\stc}{stc}
\newcommand\Crown{\textup{Cr}_4^{01}}
\newcommand\Nplu{\mathbb N^+}
\newcommand \Nnul {{\mathbb N_0}}
\newcommand{\set}[1]{\{#1\}}
\newcommand{\restrict}[2]{{#1}\mathord{\upharpoonright}_{#2}}
\newcommand\quot[1]{``#1''}
\newcommand\ivc[1]{\vec x^{(#1)}}

\begin{document}

\title{Coordinatizing coalition lattices}
\author[G.\ Cz\'edli]{G\'abor Cz\'edli}
\email{czedli@math.u-szeged.hu}
\urladdr{http://www.math.u-szeged.hu/~czedli/}
\address{University of Szeged, Bolyai Institute. 
Szeged, Aradi v\'ertan\'uk tere 1, HUNGARY 6720}

\author[E.~K.~Horv\'ath]{Eszter K.~Horv\'ath}
\email{horeszt@math.u-szeged.hu}
\urladdr{http://www.math.u-szeged.hu/~horvath/}
\address{University of Szeged, Bolyai Institute. 
Szeged, Aradi v\'ertan\'uk tere 1, HUNGARY 6720}

\date{\presentday}

\begin{abstract}
Let $P$ denote a finite poset in which no two incomparable elements have an upper bound. 
The coalition lattice $\Coal(P)$ of $P$ consists of all
subsets of $P$, ordered as follows: $X$ is less than or
equal to $Y$ if and only if there exists an injective map
from $X$ to $Y$ that sends each element to an element
greater than or equal to itself.
The main result coordinatizes $\Coal(P)$ with integer vectors. With the help of this coordinatization, the paper gives new and shorter proofs of several results on coalition lattices proved by G. Cz\'edli, Gy.\ Poll\'ak, and B.\ Larose about three decades ago. As the main application of the coordinatization, the paper proves that if $K$ is a finite simple lattice with at least three elements, then $\Coal(P)$ has no sublattice isomorphic to $K$. 
Previously, this was known only for $K=M_3$, by a result from 1996.
\end{abstract}

\keywords{Coalition order; finite poset; upper-bound-free poset; M\"obius inversion; Hall's marriage theorem; Hall--K\H{o}nig matching theorem.}

\subjclass{Primary 06B05, Secondary 06A07.}

\dedicatory{We dedicate this paper to Andrea Bozs\'o and Csilla Kert\'esz, participants in the OTKA research project No.~126581 (2017--2019), with gratitude for their invaluable support in many practical aspects of academic life before, during, and after the project.}

\maketitle

\section{Introduction}

Motivated by coalition formation in game theory and decision making, Cz\'edli and Poll\'ak \cite{CzedliPollak} introduced the concept of \emph{coalition lattices}\footnote{The concepts needed in the present paper will be defined later.} as a way of extending a partial order on a finite set \(P\) to an order on its subsets, viewed as coalitions. The underlying idea is that relative strengths of individual members can be used to compare coalitions, even in the absence of numerical measures of strength. In addition to \cite{CzedliPollak}, coalition lattices have also been studied by Cz\'edli \cite{CzedliHorn,CzedliDoubling},  by Cz\'edli, Larose, and Poll\'ak \cite{CzedliLarosePollak}, and by Davidson and Gr\"atzer \cite{davidsongratzer}. 
Although these papers establish several properties of coalition
lattices, the description of joins and meets is still complicated, and some open questions
remain unsolved. In particular, we know too little about the sublattices of coalition lattices. 

Coordinatization has always played a central role in lattice theory. For examples, see von Neumann \cite {neumann}, J\'onsson \cite{jonsson}, Day and Pickering \cite{DayPickering}, and Wehrung \cite{wehrung}. 

Our main goal is to coordinatize coalition lattices. Using the coordinatization technique that we  develop, we reprove some known properties of coalition lattices. Furthermore, by means of coordinatization, we increase our knowledge of sublattices of coalition lattices: namely, we prove that if $L$ is a finite simple lattice with at least three elements, then $L$ cannot be embedded into a coalition lattice.

Our approach was motivated by the use of filters in the papers of Horv\'ath, \v{S}e\v{s}elja, and Tepav\v{c}evi\'c \cite{HST2013,HST2015,HST2016}. Note that filters and their duals frequently occur in structural investigations in lattice theory; see, for instance, the representation theory of finite distributive lattices in Gr\"atzer \cite{GGfound}. 
The crucial advantage of the filter approach here is that it allows us to apply the Hall--K\H{o}nig matching theorem.
(We prefer the name \emph{Hall--K\H{o}nig matching theorem} to \emph{Hall's marriage theorem}, since D\'enes K\H{o}nig's equivalent theorem preceded Philip Hall's theorem by four years; see Cameron \cite{cameron} for historical details.)

\section{Coalition lattices}
Let $P=(P,\leq)$ be a partially ordered set
(or, briefly, a \emph{poset}). By convention, a poset cannot be empty. 
For subsets $X$ and $Y$ of a poset $P$, a map $f\colon X\to Y$ is said to be \emph{extensive} if $f$ is injective and $x\leq f(x)$ for all $x\in X$. Let us emphasize: an extensive map is necessarily injective (in other words, one-to-one) by definition.
By a \emph{coalition} of $P$ we mean a subset of $P$. Let $\Coal(P)$ denote the set of all coalitions of $P$. Thus, as a set, $\Coal(P)$ is the power set of $P$, but we equip it with a partial order different from set inclusion. Namely, for $X,Y\in\Coal(P)$, we define
\begin{equation}
X\cleq Y \iff \text{there exists an extensive map }f\colon X\to Y.
\label{eq:cleq}
\end{equation}

Our convention throughout is that \emph{every poset} occurring in the paper is assumed to be \emph{finite}, regardless of whether finiteness is mentioned again for emphasis or not.
Hence it is straightforward to see, and it is noted in Cz\'edli and Poll\'ak \cite{CzedliPollak}, that \eqref{eq:cleq} defines a partial order. We always equip $\Coal(P)$ with \emph{this} partial order, and from now on, we consider $\Coal(P)=\bigl(\Coal(P),\cleq\bigr)$ a poset. 

There is another approach to the relation defined in \eqref{eq:cleq}. Since the intersection of any nonempty collection of partial orders on a given set is a partial order again, it is not hard to see that $\cleq$ is the smallest partial order $\rho\subseteq \Coal(P)\times \Coal(P)$ such that, for every $X,X_1,X_2,Y,Y_1,Y_2\in \Coal(P)$ and $u,v\in P$, 
\begin{enumerate}
\item[\textup{(A1)}]if $u\leq v$, then $(\set u,\set v)\in\rho$,
\item[\textup{(A2)}] if $X\subseteq Y$, then $(X,Y)\in\rho$, and
\item[\textup{(A3)}] if $X_1\cap X_2=\emptyset$, $Y_1\cap Y_2=\emptyset$, $(X_1,Y_1)\in\rho$, and $(X_2,Y_2)\in\rho$, then we have $(X_1\cup X_2, Y_1\cup Y_2)\in\rho$.
\end{enumerate}

A poset $P$ is called \emph{upper-bound-free} if whenever $a,b,c\in P$ satisfy $a\leq c$ and $b\leq c$, the elements $a$ and $b$ are comparable. 
For example, disregarding the numbers in its vertices, Figure~\ref{fig:1} depicts the Hasse diagram of an upper-bound-free poset.
Cz\'edli and Poll\'ak \cite{CzedliPollak} proved that for a finite poset $P$,
\begin{equation}
\Coal(P)\text{ is a lattice }\quad\iff\quad P\text{ is upper-bound-free}.
\label{eq:whenCoalLat}
\end{equation}
Therefore, whenever we say that $\Coal(P)$ is a coalition lattice, we know that $P$ is a finite
upper-bound-free poset. 

The covering graph of $P$ allows us to speak of
the \emph{connected components} of $P$. 
For example, the poset in Figure~\ref{fig:1} has exactly two connected components.
Observe that, for a finite upper-bound-free poset $P$, 
each connected component of $P$ has a unique smallest element.

For an element $u$ of a finite (not necessarily upper-bound-free) poset $P$, let $\UpCov(u)$ denote the \emph{set of upper covers} of $u$.
For $u\in P$, the smallest \emph{principal filter} containing $u$ is denoted by $\fil(u)=\fil_P(u)$; in notation, $\fil(u):=\set{x\in P: u\leq x}$. Dually, $\idl(u) =\idl_P(u):= \set{x\in P: x\leq u}$ is the \emph{principal ideal} generated by $u$.
The \emph{set of maximal elements} and the \emph{set of minimal elements} of $P$ are denoted by $\Max(P)$, and $\Min(P)$, respectively. That is, 
\begin{equation}
\Max(P)=\set{x\in P: |\fil(x)|=1} \text{ and } \Min(P)=\set{x\in P: |\idl(x)|=1}.
\label{eq:MaxMinP}
\end{equation}

\section{Coordinatization}

In this section, $P=(P,\leq)$ denotes a finite upper-bound-free poset.

By a \emph{coordinatization} of $\Coal(P)$, we mean an isomorphism from $\Coal(P)$ onto a structure whose elements are tuples with entries in a ring. We call the latter structure the \emph{coordinate space} of $P$. Since it consists of the so-called eligible tuples (to be defined later), we denote the coordinate space of $P$ by $\ET(P)$. Its definition needs some preparation.

For $n\in\Nplu:=\set{1,2,3,4,\dots}$, we denote the set $\set{1,2,\dots,n}$ by $[n]$. 
An $n$-tuple $(x_1,\dots,x_n)$ is just another notation for a map from $[n]$ to a set containing $x_1,\dots,x_n$.  It will be convenient to denote this $n$-tuple also in another way:
\begin{equation}
\text{let }\quad(x_i:i\in[n]) \quad \text{ mean }\quad (x_1,x_2,\dots,x_n)=:\vec x.
\label{eq:ntuple}
\end{equation}

Conversely, any map with domain $[n]$ can be considered an $n$-tuple. 
We extend this tuple notation to maps whose domain is $P$: let $f\colon P\to A$ be a map.
Then, motivated by \eqref{eq:ntuple} and 
without explicitly specifying the codomain $A$ of $f$,
\begin{equation}
\text{we may use }(f(q):q\in P)\text{ to denote the map }f\colon P\to A.
\label{eq:Ptuple}
\end{equation}
Furthermore, we can visualize $(f(q):q\in P)$ by labeling the vertices $q\in P$ in the diagram of $P$ with the corresponding values $f(q)$.

For $q\in P$, since $P$ is finite and upper-bound-free,
we may compute the \emph{size}
$\spf(q)=\spf_P(q)$ of the \emph{principal filter} $\fil(q)$
generated by $q$ recursively as follows:
\begin{equation}
\spf(q)=\spf_P(q):=\begin{cases}
1,&\text{if }q\in\Max(P),\cr
1+\sum_{b\in\UpCov(q)}\spf(b),&\text{otherwise}.
\end{cases}
\label{eq:filternumber}
\end{equation}
Indeed, using \eqref{eq:MaxMinP}, a straightforward downward induction on $P$ yields that
\begin{equation}
\spf(q)=|\fil(q)|\qquad\text{for all }q\in P.
\label{eq:fildef}
\end{equation}

\begin{figure}[ht] 
\centerline{ \includegraphics[width=\figwidthcoeff\textwidth]{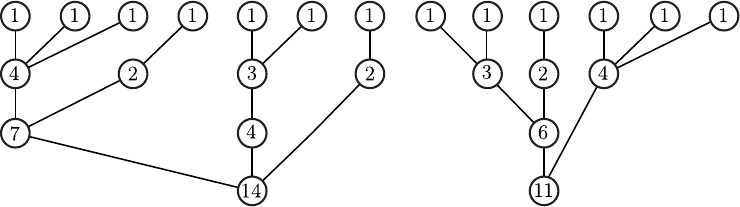}} \caption{An upper-bound-free poset $P$ and the $P$-tuple $\spf_P$}\label{fig:1}
\end{figure}

According to \eqref{eq:Ptuple},  $\spf_P$ is a $P$-tuple. In Figure \ref{fig:1}, 
the circles stand for the elements of $P$, while the numbers in these circles represent 
the $P$-tuple $\spf$. 
As usual, an empty sum is understood to be \(0\) in the following definition.

\begin{definition}\label{def:elig}
We say that a $P$-tuple $\vec x=(x_q: q\in P)$ of nonnegative integers is \emph{eligible} (with respect to $P$) if for each $q\in P$,
\begin{equation}
 \sum_{b\in\UpCov(q)} x_b  \leq x_q \leq 1+\sum_{b\in\UpCov(q)} x_b.
\label{eq:elIg}
\end{equation}
The \emph{set of eligible $P$-tuples} of nonnegative integers is denoted by $\ET(P)$. We equip it with the coordinate-wise order $\eleq$: for $\vec x=(x_q:q\in P)$  and $\vec y=(y_q:q\in P)$ in $\ET(P)$, we let
\begin{equation*}
\vec x\eleq \vec y \iff x_q\leq y_q\text{ for all }q\in P.
\end{equation*}
The \emph{coordinate space} associated with $P$ is the poset $\ET(P)=\bigl(\ET(P),\eleq)$.
\end{definition}
Comparing \eqref{eq:filternumber} and \eqref{eq:elIg}, we immediately obtain that
\begin{equation}
(0: q\in P)=0_{\ET(P)}\text{ \ and \ } \spf_P=1_{\ET(P)},
\label{eq:ETbounded}
\end{equation}
and therefore $\ET(P)$ is a bounded poset. Another immediate consequence of \eqref{eq:ETbounded} is that for every $\vec x\in\ET(P)$ and $q\in P$, we have $0\leq x_q\leq \spf_P(q)$.
It will soon appear that  $\ET(P)$ is a lattice.
Before formulating our main result, we still need another definition.

\begin{definition}\label{def:reciprmaps} Let $P$ be a finite upper-bound-free poset. We define two maps, the \emph{coordinatizing map} $\coord$ and the \emph{reconstruction map} $\rec$ as follows.
\begin{align}
&\text{for }X\in \Coal(P),\text{ let }\coord(X):=\bigl(|X\cap\fil(q)|: q\in P \bigr),\label{eq:coordmap}\\
&\text{for }\vec x=(x_q:q\in P)\in\ET(P),\text{ let }\rec(\vec x):=\set{q: x_q=1+\sum_{b\in\UpCov(q)}x_b}.
\label{eq:recmap}
\end{align}
For $X\in\Coal(P)$, we call $\coord(X)$ the \emph{coordinate tuple} (or the \emph{coordinate vector}) of $X$. We may subscript $\coord$ and $\rec$ with $P$ if ambiguity threatens.
\end{definition}

\begin{figure}[ht] 
\centerline{ \includegraphics[width=\figwidthcoeff\textwidth]{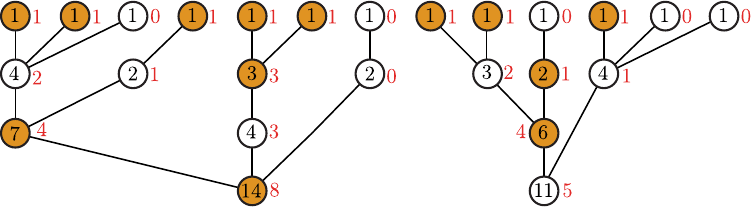}} \caption{The red numbers form an eligible $P$-tuple }\label{fig:2}
\end{figure}

In \eqref{eq:recmap}, the equality $x_q=1+\sum_{b\in\UpCov(q)}x_b$ can equivalently be replaced by $x_q>\sum_{b\in\UpCov(q)}x_b$ or $x_q\ne\sum_{b\in\UpCov(q)}x_b$. For the 25-element poset $P$ drawn in Figure \ref{fig:1}, Definition \ref{def:reciprmaps} is visualized by Figure \ref{fig:2}. 
 Namely, if $X$ consists of the elements represented by the orange-filled circles, then $\coord(X)$ is the system of the red numbers (that is, numbers outside the little circles). On the other hand, if $\vec x$ is the system of the red numbers, then $\rec(\vec x)$ consists of the elements represented by the orange-filled circles.

We are now in a position to state the first theorem of the paper.

\begin{theorem}\label{thm:isomorphic} Let $P$ be a finite upper-bound-free poset. Then 
\[\coord\colon\Coal(P)\to\ET(P)\quad\text{ and } 
\quad \rec\colon\ET(P)\to\Coal(P)
\]
are mutually inverse order isomorphisms.
\end{theorem}

\begin{remark}
By a result of Cz\'edli and Poll\'ak \cite{CzedliPollak}, $\Coal(P)$ above is a lattice. However, this fact is not used 
in Theorem \ref{thm:isomorphic} or its proof, since we will use the theorem to give a simple new proof of the old result.
\end{remark}

\begin{remark}
For a large finite upper-bound-free poset $P$, deciding whether two members of $\Coal(P)$ are comparable directly from \eqref{eq:cleq} can be computationally demanding. In contrast, 
deciding whether $\vec x\leq_{\mathrm{et}}\vec y$ for two members $\vec x,\vec y\in\ET(P)$ requires at most $|P|$ comparisons of integers. This computational advantage alone provides a motivation for Theorem~\ref{thm:isomorphic}.
\end{remark}

For a map $f\colon A\to B$ and a subset $C\subseteq A$, we denote by $\restrict f C$ the \emph{restriction} of $f$ to $C$.

\begin{proof}[Proof of Theorem \ref{thm:isomorphic}]
Assume that $q\in P$ and $X\in\Coal(P)$. Let $\vec x=(x_r:r\in P):=\coord(X)$. Since the elements in $\UpCov(q)$ are pairwise incomparable and  $P$ is finite and upper-bound-free, we obtain that  
\begin{align}
&\fil(q)=\set q \cup \bigcup_{b\in\UpCov(q)}\fil(b)\quad\text{ and}
\label{eq:filquniona}\\
&\fil(b_1)\cap\fil(b_2)=\emptyset\text{ for distinct }b_1,b_2\in\UpCov(q). 
\label{eq:filqunionb}
\end{align}
Note that for $q\in\Max(P)$, $\UpCov(q)=\emptyset$ and the union in \eqref{eq:filquniona} is $\emptyset$.
Since intersection distributes over union and preserves disjointness, \eqref{eq:filquniona} and \eqref{eq:filqunionb} yield that
\begin{align}
x_q&= |X\cap\fil(q)|=|X\cap\set{q}|+\bigl| \bigcup_{b\in\UpCov(q)}(X\cap\fil(b))\bigr| 
\label{eq:cmpxqa}\\
&=|X\cap\set{q}|+ \sum_{b\in\UpCov(q)}|X\cap\fil(b)|=|X\cap\set{q}|+ \sum_{b\in\UpCov(q)}x_b. \label{eq:cmpxqb}
\end{align}
Since $|X\cap\set{q}|\in\set{0,1}$, \eqref{eq:cmpxqa} and \eqref{eq:cmpxqb} show that $\vec x$ satisfies \eqref{eq:elIg}, and thereby  $\coord(X)=\vec x\in\ET(P)$. Hence $\coord$ is indeed a $\Coal(P)\to\ET(P)$ map.  Trivially, $\rec$ is an $\ET(P)\to\Coal(P)$ map. 

Let $X\in\Coal(P)$ and  $q\in P$, and denote $\coord(X)$ by $\vec x$. Using  \eqref{eq:coordmap}, \eqref{eq:recmap} \eqref{eq:filquniona}, and \eqref{eq:cmpxqa}--\eqref{eq:cmpxqb}, 
we obtain that
\begin{align*}
q\in X &\iff |X\cap\fil(q)|=1+ \bigl|\bigcup_{b\in\UpCov(q)}(X\cap\fil(b))\bigr| 
\cr
&\iff x_q=1+\sum_{b\in\UpCov(q)}x_b\iff q\in\rec(\vec x)=\rec(\coord(X)).
\end{align*}
Hence $\rec(\coord(X))=X$, and so the composite map $\rec\circ \coord$ is the identity map of $\Coal(P)$.

Proceeding downwards, we define the \emph{layers} $H_0, H_1, \dots, H_t$ of $P$ by induction. Let $H_0:=\Max(P)$. 
Then for $i\in\Nnul:=\set{0,1,2,3,\dots}$, we let 
\[H_{i+1}:=\Max(P\setminus(H_0\cup\dots\cup H_i))
\] 
(with the convention that $\Max(\emptyset)=\emptyset$), and let $t$ be the largest $i\in\Nnul$ such that $H_i\neq\emptyset$. Then $\set{H_0,\dots,H_t}$ is a partition of $P$.  For an element $s\in P$, the \emph{depth} of $s$  is the unique subscript  $i$ such that $s\in H_i$.

Next, let $\vec y\in\ET(P)$.  Denote $\rec(\vec y)$ by $Y$ and $\coord(Y)$ by $\vec z$.   Let $q\in P$; by induction on $\dep(q)$, we are going to show that $y_q=z_q$. If $\dep(q)=0$, that is, if $q\in\Max(P)$, then $y_q\in \set{0,1}$, since $\vec y$ is eligible, and $z_q\in\set{0,1}$ by \eqref{eq:coordmap}. Furthermore, $y_q=1\iff q\in Y\iff z_q=1$, and so $y_q=z_q$. Now let $\dep(q)>0$, and assume that $y_b=z_b$ for all $b$ with $\dep(b)<\dep (q)$; in particular, for all  $b\in\UpCov(q)$.  Hence, letting 
\[s_q(\vec y):=\sum_{b\in\UpCov(q)} y_b\quad\text{ and }\quad s_q(\vec z):=\sum_{b\in\UpCov(q)} z_b,
\]
$s_q(\vec y)=s_q(\vec z)$. By the eligibility of $\vec y$, we have $y_q\in\set{s_q(\vec y), 1+s_q(\vec y)}$.  If $y_q=1+s_q(\vec y)$, then $q\in Y=\rec(\vec y)$ by \eqref{eq:recmap}, and $z_q=1+s_q(\vec z)=1+s_q(\vec y)=y_q$ by \eqref{eq:coordmap} and \eqref{eq:filquniona}. Similarly, if 
$y_q=s_q(\vec y)$, then $q\notin Y=\rec(\vec y)$ by \eqref{eq:recmap}, and $z_q=s_q(\vec z)=s_q(\vec y)=y_q$ by \eqref{eq:coordmap} and \eqref{eq:filquniona}. Hence $\vec y=\vec z$, and we have seen that the composite map $\coord\circ \rec$ is the identity map of $\ET(P)$. Therefore, $\coord$ and $\rec$ are mutually inverse bijections, and we need to show that both are order-preserving.

First, assume that $X,Y\in\Coal(P)$ such that $X\cleq Y$. Denote $\coord(X)$ and $\coord Y$ by $\vec x$ and $\vec y$, respectively.
Let $f\colon X\to Y$ be an extensive map, and consider an arbitrary $q\in P$.  For $x\in X\cap \fil(q)$,  we have $f(x)\in Y\cap\fil(q)$, since $f$ is extensive. Therefore, since $f$ is injective, so is the restriction 
$\restrict f{X\cap \fil(q)}\colon X\cap \fil(q) \to Y\cap\fil(q)$. Consequently, 
\[x_q=|X\cap \fil(q)| \leq |Y\cap \fil(q)| = y_q.
\]
Since $q$ was an arbitrary element of $P$, we have $\coord(X)=\vec x\eleq \vec y=\coord(Y)$. Therefore, $\coord$ is order preserving.

Second, assume that $\vec x,\vec y\in\ET(P)$ such that $\vec x\eleq \vec y$. Let $X:=\rec(\vec x)$ and $Y:=\rec(\vec y)$. Applying $\coord$, which is the inverse of $\rec$, to both sides of these equalities, we obtain that
$\coord(X)=\vec x$ and $\coord(Y)=\vec y$.
Consider the bipartite graph $G$ with partite sets $X\times\set{1}$ and $Y\times\set{2}$ and edge set
\begin{equation}
\bigl\{ \set{(x,1), (y,2)}: x\in X,\ y\in Y,\text{ and }x\leq y\bigr\}.
\label{eq:bipartite}
\end{equation}
(We use $\times \set{1}$ and $\times \set{2}$ above to ensure that the partite sets are disjoint.)
For $S\subseteq X$, let 
\begin{equation}
N(S):=\set{y\in Y: \text{there exists }x\in S\text{ such that }x\leq y}. 
\label{eq:NSdef}
\end{equation}
Clearly, $|N(S)|$ is the number of neighbours of $S\times\set{1}$ in the bipartite graph $G$; see \eqref{eq:bipartite}. By finiteness and transitivity, we can reformulate \eqref{eq:NSdef} to obtain that
\begin{equation*}
N(S):=\set{y\in Y: \text{there exists }q\in \Min(S)\text{ such that } y\in\fil_P(q)}. 
\end{equation*}
Formulating this in a more convenient way, 
\begin{equation}
N(S):=\bigcup_{q\in\Min(S)}( Y\cap\fil_P(q)).
\label{eq:MinNSdef}
\end{equation}
If $q_1$ and $q_2$ are distinct elements of $\Min(S)$, then they are incomparable, and  therefore $\fil_P(q_1)$ is disjoint from $\fil_P(q_2)$, since $P$ is upper-bound-free. 
Hence the union in \eqref{eq:MinNSdef} is the union of disjoint sets. In the computation below we use, in this order, \eqref{eq:MinNSdef} together with the just-mentioned disjointness, \eqref{eq:coordmap},  $\vec x\eleq \vec y$, \eqref{eq:coordmap} again, 
$X\supseteq S$, the disjointness again, and the obvious fact that each element of $S$ belongs to $\fil_P(q)$ for some $q\in\Min (S)$:
\begin{align*}
|N(S)|&=\sum_{q\in\Min(S)}| Y\cap\fil_P(q)|=\sum_{q\in\Min(S)} y_q \geq \sum_{q\in\Min(S)} x_q 
\cr
&=\sum_{q\in\Min(S)}| X\cap\fil_P(q)| \geq \sum_{q\in\Min(S)}| S\cap\fil_P(q)| 
\cr
 &=\bigl|\bigcup_{q\in\Min(S)} (S\cap\fil_P(q)) \bigr|
 =\bigl|S\cap \bigcup_{q\in\Min(S)} \fil_P(q) \bigr| = |S|.
 \end{align*}
Hence the number of neighbours of $S\times\set1$ in the graph $G$ is at least $|S|=|S\times \set 1$.  

Since $S$ was an arbitrary subset of $X$, we have seen that an arbitrary subset of the partite set $X\times \set 1$ has at least as many neighbours in $Y\times\set{2}$ as the size of the subset in question.  Hence the Hall--K\H{o}nig matching theorem yields an injective map
$f\colon X\times\set{1}\to Y\times\set 2$ such that, 
for each $(x,1)\in X\times\set 1$,
$\{(x,1), f((x,1))\}$ is an edge of $G$. Define a map $g\colon X\to Y$ by the rule $y=g(x)\iff (y,2)=f((x,1))$. 
Then $g$ is injective and, by \eqref{eq:bipartite},  $x\leq g(x)$ for all $x\in X$. Hence $g\colon X\to Y$ is an extensive map, and $\rec(\vec x)=X\cleq Y=\rec(\vec y)$. This shows that, similarly to $\coord$,  $\rec$ is also order-preserving. Hence $\coord$ and $\rec$ are inverse order isomorphisms, completing the proof of Theorem \ref{thm:isomorphic}.
\end{proof}

\section{The coordinate space is a lattice}

The implication "$\Longrightarrow$" of \eqref{eq:whenCoalLat} is easy to prove, whereas the converse implication is less straightforward. In addition to the original proof by Cz\'edli and Poll\'ak \cite{CzedliPollak}, two different proofs of the latter implication were given by Cz\'edli, Larose, and Poll\'ak \cite{CzedliLarosePollak} and by Davidson and Gr\"atzer \cite{davidsongratzer}. The following theorem is a consequence of \eqref{eq:whenCoalLat} and Theorem~\ref{thm:isomorphic}. However, we give a direct proof of it, thereby obtaining a fourth proof of the "$\Longleftarrow$" part of \eqref{eq:whenCoalLat}.

\begin{theorem}\label{thm:ETlat}
If $P$ is a finite upper-bound-free poset, then $\ET(P)$ is a lattice. For $\vec x,\vec y\in\ET(P)$, their meet $\vec u:=\vec x\wedge\vec y$ and join $\vec v:=\vec x\vee\vec y$ are given recursively, proceeding from the maximal elements of $P$ downwards, by
\begin{align}
\vec u&=\bigl(\min\{x_q,y_q,1+\sum_{b\in\UpCov(q)}u_b\}:q\in P\bigr),
\label{eq:ETPmeet}\\
\vec v&=\bigl(\max\{x_q,y_q,\sum_{b\in\UpCov(q)}v_b\}:q\in P\bigr).
\label{eq:ETPjoin}
\end{align}
\end{theorem}

\begin{proof} First, we focus on $\vec u$. 
For $r\in P$, denote $\sum_{b\in\UpCov(r)}u_b$  by $S_r$. 
Note the rule: $u_r=\min\set{x_r, y_r, 1+S_r}$, for every $r\in P$. 

Let $q\in P$. 
Since $u_b\leq x_b$ for all $b\in P$ by \eqref{eq:ETPmeet}, \eqref{eq:elIg} implies that 
$S_q\leq \sum_{b\in\UpCov(q)}x_b\leq x_q$.
We obtain $S_q\leq y_q$ similarly. Hence 
\[
S_q \leq \min\set{x_q,y_q,1+S_q} = u_q\quad\text{ and, clearly, }\quad u_q \leq 1+S_q.
\]
By \eqref{eq:elIg}, since $q\in P$ was arbitrary, this means that $\vec u\in \ET(P)$.
By \eqref{eq:ETPmeet}, $\vec u\eleq \vec x$ and $\vec u\eleq \vec y$. 

Assume now that $\vec z\in \ET(P)$ is a common lower bound of $\vec x$ and $\vec y$, that is, $\vec z\eleq \vec x$ and $\vec z\eleq \vec y$. We show by induction on $\dep(r)$ that 
\begin{equation}
\text{for every \ }r\in P,\text{ \ we have \ } z_r\leq u_r.
\label{eq:uzrind}
\end{equation}
Let $\dep(r)=0$, that is, $r\in\Max(P)$. Then, by \eqref{eq:elIg} and since the empty sum is $0$, 
\begin{equation}
x_r\leq 1+\sum_{b\in\UpCov(r)}x_b=1\leq 1+S_r
\text{ \ and, similarly, \ }y_r\leq 1+S_r.
\label{eq:xqyqSq}
\end{equation}
Combining $z_r\leq x_r$, $z_r\leq y_r$, and \eqref{eq:xqyqSq}, we obtain 
\[
z_r \leq\min\set{x_r,y_r} = \min\set{x_r,y_r,1+S_r}=u_r,
\]
as required. 

Next, assume that $\dep(r)>0$ and $z_b\leq u_b$ for all $b\in P$ with $\dep(b)<\dep(r)$. Then
$\sum_{b\in\UpCov(r)} z_b\leq S_r$. Hence, by \eqref{eq:elIg}, $z_r\leq 1+\sum_{b\in\UpCov(r)} z_b\leq 1+S_r$. Combining this inequality with $z_r\leq x_r$ and $z_r\leq y_r$, we obtain that 
$z_r\leq \min\set{x_r,y_r,1+S_r}=u_r$,
and we have verified \eqref{eq:uzrind}.
Therefore, $\vec z\eleq \vec u$, and so $\vec u$ is the greatest lower bound of $\vec x$ and $\vec y$. This proves that $\vec u$ is the meet of $\vec x$ and $\vec y$.

Note that we have already proved that $\ET(P)$ is a lattice, since it has a largest element, $(\spf(q): q\in P)$, and every finite meet-semilattice with a largest element is a lattice.

Next, we deal with $\vec v$. For $r\in P$, we denote $\sum_{b\in\UpCov(r)}v_b$ by $T_r$. Then our computational rule is $v_r=\max\set{x_r,y_r,T_r}$.  

Let $q\in P$.
By \eqref{eq:elIg},  $x_q\leq 1+\sum_{b\in\UpCov(q)} x_b$.
Since $x_b\leq v_b$ by \eqref{eq:ETPjoin}, the previous inequality implies $x_q\leq 1+T_q$. We obtain $y_q\leq 1+T_q$ similarly. Hence $v_q=\max\set{x_q,y_q,T_q}\leq 1+T_q$, and we also have $T_q\leq v_q$. Thus, by \eqref{eq:elIg}, $\vec v\in\ET(P)$. 
By \eqref{eq:ETPjoin}, $\vec x\eleq \vec v$ and $\vec y\eleq \vec v$.

Assume that $\vec w\in \ET(P)$ satisfies $\vec x\eleq \vec w$ and $\vec y\eleq \vec w$; we need to show that $\vec v\eleq \vec w$. We show by induction on $\dep(q)$ that
\begin{equation}
\text{for all \ }q\in P,\text{ \ we have \ }v_q\leq w_q.
\label{eq:indvrwr}
\end{equation}
If $\dep(q)=0$ then $T_q$ is the empty sum, and therefore
\[v_q=\max\set{x_q, y_q,T_q}=\max\set{x_q, y_q,0}\leq w_q,
\]
as required. 
Assume that $\dep(q)>0$ and the inequality in \eqref{eq:indvrwr}, with $b$ instead of $q$, holds for all $b\in P$ with $\dep(b)<\dep(q)$. Then, using $v_b\leq w_b$ for $b\in\UpCov(q)$ and \eqref{eq:elIg} for $\vec w$, we have that
\begin{align*}
v_q&=\max\set{x_q,y_q,T_q}\leq \max{\set{w_q,w_q, \sum_{b\in \UpCov(q)}w_b }} =w_q.
\end{align*}
This proves \eqref{eq:indvrwr}, and so $\vec v\eleq\vec w$.
Thus $\vec v=\vec x\vee \vec y$, and the proof of  Theorem \ref{thm:ETlat} is complete.
\end{proof}

\section{New proofs of some known results}

Once Theorem \ref{thm:ETlat} (or at least its  \eqref{eq:ETPmeet} part) and Theorem~ \ref{thm:isomorphic} have been established, the proofs in this section are shorter than their original counterparts.

We recall the following fact from Cz\'edli and Poll\'ak \cite{CzedliPollak}, which is straightforward to prove: for a finite poset $P$,
\begin{equation}
\text{if }P\text{ is not upper-bound-free, then }\Coal(P)
\text{ is not a lattice.}
\label{eq:notUBFnotlat}
\end{equation}
Indeed---slightly modifying the argument given in Cz\'edli and Poll\'ak \cite{CzedliPollak}---suppose for a contradiction that $\Coal(P)$ is a lattice but $P$ is not upper-bound-free. Finiteness allows us to pick $a,b,c\in P$ such that $a\neq b$, $a \prec c$, and $b\prec c$.
(Indeed, a minimal upper bound of two incomparable elements has at least two lower covers.) Let $X:=\set a\vee\set b$ in $\Coal(P)$. Since $\set{a}\cleq \set{c}$ and $\set b\cleq \set c$, we have $\set a\cleq X\cleq \set c$ and $\set b\cleq X\cleq \set c$. Hence \eqref{eq:cleq} yields that $X$ is a singleton. Thus write $X=\set x$, and we have $a\leq x\leq c$ and $b\leq x\leq c$. Therefore, since $a\prec c$ and $b\prec c$, we have $x=c$, that is, $X=\set c$. 
Since $\set a\cleq \set{a,b}$ and $\set b\cleq \set{a,b}$, we obtain that 
$\set c=\set a\vee \set b\cleq \set{a,b}$. Hence it follows from \eqref{eq:cleq} that $c\leq a$ or $c\leq b$.
This is a contradiction proving \eqref{eq:notUBFnotlat}.  

Now we formulate the main result of Cz\'edli and Poll\'ak \cite{CzedliPollak} in a slightly less general but simpler form than the original statement.

\begin{corollary}[Cz\'edli and Poll\'ak \cite{CzedliPollak}]
Let $P$ be a finite poset. Then $\Coal(P)$ is a lattice if and only if $P$ is upper-bound-free.
\end{corollary}

\begin{proof} The \quot{only if} implication is \eqref{eq:notUBFnotlat}. The \quot{if} implication follows from Theorems \ref{thm:isomorphic} and \ref{thm:ETlat}.
\end{proof}

We know from Cz\'edli and Poll\'ak \cite[Proposition 2]{CzedliPollak} that for every finite poset $P$  (and even for every finite quasiordered set, which we do not define here), the set-theoretic complementation is a dual order isomorphism of $\Coal(P)$. Below we deal only with a particular case of this result.

\begin{corollary}[Cz\'edli and Poll\'ak \cite{CzedliPollak}]
\label{corol:dual}
Let $P$ be a finite upper-bound-free poset. Then the map
\begin{equation*}
\stc\colon \Coal(P)\to\Coal(P),\text{ defined by }
X\mapsto P\setminus X
\end{equation*}
is a dual automorphism of $\Coal(P)$. Hence $\Coal(P)$ is a selfdual lattice.
\end{corollary}

The acronym $\stc$ comes from \quot{set-theoretic complementation}. Note that, for $X\in \Coal(P)$, $\stc(X)$ need not be a lattice-theoretic complement of $X$ in $\Coal(P)$. For example, if
$P=\set{u,v}$ is the two-element chain with $u\prec v$, then
$\stc(\set u)=\set v$, and so  $\set u\vee \stc(\set u)=\set v$ is different from $1_{\Coal(P)}=P$.

\begin{proof}[Proof of Corollary \ref{corol:dual}] Clearly, $\stc\colon \Coal(P)\to\Coal(P)$ is an involutive bijection.
For $\vec x,\vec y\in\ET(P)$, addition and subtraction are understood componentwise. 
It is clear by \eqref{eq:fildef} and 
\eqref{eq:coordmap} that, for $X\in \Coal(P)$,
$\coord(X)+\coord(\stc(X))=\spf$, that is,
\begin{equation}
\coord(\stc(X))=\spf-\coord(X).
\label{eq:coordfilcoord}
\end{equation} 
In particular, it follows that $\spf-\coord(X)\in\ET(P)$
for every $X\in\Coal(P)$. 
Hence, we could use $\eleq$ instead of the componentwise $\leq$ in the third comparison of the display below.  
(We use $\leq$, however, to suggest that subtraction is somewhat alien to the lattice $\ET(P)$.)
For $X,Y\in\Coal(P)$, Theorem \ref{thm:isomorphic} and \eqref{eq:coordfilcoord} yield that
\begin{align*}
X\cleq Y&\iff \coord(X)\eleq\coord(Y) 
\cr
&\iff \spf-\coord(Y)\leq \spf-\coord(X)
\cr 
&\iff \coord(\stc(Y)) \eleq \coord(\stc(X))
\iff \stc(Y)\cleq \stc(X),
\end{align*}
which completes the proof.
\end{proof}

Let $P$ be a finite upper-bound-free poset. We say that $X\in \Coal(P)$ is a \emph{winning coalition} if $\stc(X)\cleq X$. The proof of the following statement, including the arguments for three required lemmas, is about one and a half pages long in Cz\'edli, Larose, and Poll\'ak \cite{CzedliLarosePollak}. Here we will give a substantially shorter proof.

\begin{corollary}[Cz\'edli, Larose, and Poll\'ak \cite{CzedliLarosePollak}]\label{corol:winning}
Let $P$ be a finite upper-bound-free poset, and let $X,Y\in\Coal(P)$. If $X$ and $Y$ are winning coalitions, then so is $X\wedge Y$. Furthermore, the winning coalitions form a principal filter in $\Coal(P)$.
\end{corollary}

\begin{proof} We will use Theorem \ref{thm:isomorphic} implicitly. For $\vec x\in\ET(P)$, let $\halfspf\leq \vec x$ mean that $\spf(q)/2\leq x_q$ for all $q\in P$.  (We cannot write $\halfspf \eleq \vec x$ since $\halfspf\notin\ET(P)$.) 

It follows from \eqref{eq:coordfilcoord} that a coalition $X\in \Coal(P)$ is winning if and only if $\halfspf\leq \coord(X)$. Using this fact, we may work only in $\ET(P)$. 
We call a $P$-tuple $\vec x\in\ET(P)$ \emph{winning} if $\halfspf\leq \vec x$. By the second half of  \eqref{eq:ETbounded}, $\vec x$ is winning if and only if $\halfspf\leq \vec x\eleq \spf$. 
Assume that $\vec x,\vec y\in \ET(P)$ are winning $P$-tuples, and let $\vec u:={\vec x}\wedge{\vec y}$; see \eqref{eq:ETPmeet}. To prove that $\vec u$ is winning, we show by induction on $\dep(q)$ that $\spf(q)/2\leq u_q$ for all $q\in P$. If $\dep(q)=0$, then $1/2=\spf(q)/2\leq x_q\leq \spf(q)=1$ and $x_q\in \Nnul$ imply $x_q=1$. Similarly, $y_q=1$.
By \eqref{eq:ETPmeet}, $u_q=1$. Hence, if $\dep(q)=0$, then $\spf(q)/2\leq u_q$, as required. 
Assume now that $\dep(q)>0$ and $\spf(b)/2\leq u_b$ for all $b\in P$ with $\dep(b)<\dep(q)$. Then, by \eqref{eq:filternumber},
\begin{equation*}
\frac{\spf(q)}2
=\frac12+\frac12\sum_{b\in\UpCov(q)}\spf(b)
\le\frac12+\sum_{b\in\UpCov(q)}u_b
\le 1+\sum_{b\in\UpCov(q)}u_b.
\end{equation*}
This inequality, $\spf(q)/2\leq x_q$, $\spf(q)/2\leq y_q$, and \eqref{eq:ETPmeet} imply $\spf(q)/2\leq u_q$. Hence $\vec u$ is winning.

Let $W:=\{\vec x\in \ET(P): \vec x$ is winning$\}$. We have seen that $W$ is meet-closed. Since $\spf\in W$, $W\neq\emptyset$. Furthermore, $W$ is clearly an order filter, that is,   if $\vec x\in W$, then $\fil_{\ET(P)}(\vec x)\subseteq W$. These facts and finiteness imply that 
$W$ is a principal filter of $\ET(P)$. Hence Theorem \ref{thm:isomorphic} completes the proof of Corollary \ref{corol:winning}.
\end{proof}

The following statement uses a notation introduced in \eqref{eq:MaxMinP}, and $\prod$ denotes the direct product.

\begin{corollary}[Cz\'edli and Poll\'ak \cite{CzedliPollak}]\label{corol:dirprod} If $P$ is a finite upper-bound-free poset, then
\begin{equation}
\Coal(P)\cong\prod_{u\in\Min(P)}\Coal({\fil_P(u)})
\text{ \ and \ }
\ET(P)\cong\prod_{u\in\Min(P)}\ET({\fil_P(u)})
\label{eq:dirprod}
\end{equation}
\end{corollary}

This corollary is straightforward, and thus no proof will be given here; we only note that the $\ET(P)$ part of \eqref{eq:dirprod} implies its $\Coal(P)$ part by Theorem \ref{thm:isomorphic}, and vice versa. The proof of the $\Coal(P)$ part (even without assuming that $P$ is upper-bound-free) is also straightforward in Cz\'edli and Poll\'ak \cite{CzedliPollak}.

Note that Cz\'edli and Poll\'ak \cite{CzedliPollak}, 
Cz\'edli, Larose, and Poll\'ak \cite{CzedliLarosePollak}, and Cz\'edli \cite{CzedliHorn,CzedliDoubling}
contain further results on $\Coal(P)$ that are not mentioned in the present paper. Some of these results could probably be reproved in the $\ET(P)$ setting, but this aim is not pursued in the present paper.

\section{A non-embeddability result}
As usual, $M_3$ denotes the five-element modular 
lattice with three atoms, which are also coatoms.
It has been proved in Cz\'edli \cite{CzedliHorn} that for every finite upper-bound-free poset $P$, $M_3$ cannot be embedded into $\Coal(P)$.
Here we prove the following more general statement. A lattice is \emph{simple} if it has at most two congruences.
For example, $M_3$ is a simple lattice.

\begin{proposition}\label{prop:simple} If $P$ is a finite upper-bound-free poset and $K$ is a finite simple lattice with at least three elements, then $K$ cannot be embedded into $\Coal(P)$.
\end{proposition}

Many known finite simple lattices $K$ have $M_3$ as a sublattice; for these $K$, Proposition~\ref{prop:simple} does not yield anything new. However, there are finite simple lattices without an $M_3$ sublattice. For example, the ten-element crown lattice $\Crown$, with atoms $a_i$ and coatoms $b_i$, $i\in\mathbb Z_4$, whose covering relation satisfies $a_i\prec b_j$ if and only if $j=i$ or $j=i+1$ in $\mathbb Z_4$, is such a lattice. This can easily be verified using rotational symmetry and, for example, Gr\"atzer's technical lemma \cite{GGtechnlemma}.

\begin{proof}[Proof of Proposition \ref{prop:simple}]
Theorem \ref{thm:isomorphic} allows us to deal with coordinate spaces $\ET(P)$ rather than with coalition lattices. We may assume that $|P|\geq 2$, since $|\ET(P)|=2$ when $|P|=1$.
Aiming towards a contradiction, suppose that $K$ 
is an at least three-element finite simple lattice and  there exists a finite upper-bound-free poset $P$ such that $K$ is isomorphic to a sublattice of $\ET(P)$. Then we may select such a $P$ of minimum size. Thus we may assume that $K$ is a sublattice of $\ET(P)$ but for any finite upper-bound-free poset $Q$,
\begin{equation}
\text{if }|Q|<|P|,\text{ then }K\text{ is not isomorphic to a sublattice of }\ET(Q).
\label{eq:Pminsize}
\end{equation}

For $u\in\Min(P)$, denote $\ET(\fil_P(u))$ by $L_u$. By \eqref{eq:dirprod}, we may assume that 
$K$ is a sublattice of $L:=\prod_{u\in\Min(P)} L_u$. Let $\pi_u\colon L\to L_u$ denote the natural projection and let $\Theta_u:=\{(x,y)\in K^2: \pi_u(x)=\pi_u(y)\}$ stand for the congruence kernel of the restricted map $\restrict{\pi_u}K$. Since $K$ is simple, $\Theta_u$ is either the equality or the universal congruence $K^2$. 

If $\Theta_u$ were the universal relation for all $u\in\Min(P)$, then $\pi_u(K)$ would be a singleton for all $u\in\Min(P)$, and so $K$ would also be a singleton, contradicting $|K|\geq 3$. Hence there exists an element $u\in\Min(P)$ such that $\Theta_u$ is the equality relation, and thus $\restrict{\pi_u}K$ is injective. Thus $K$ is embeddable into the direct factor $L_u=\ET(\fil_P(u))$. By \eqref{eq:Pminsize}, this is possible only if $\fil_P(u)=P$. Consequently, $u$ is the unique minimal element of $P$; denote it by $0_P$.
Furthermore, instead of being embedded, we may assume that $K$ is a sublattice of $\ET(P)$.

Let $t:=|K|$, let $K=\set{\ivc 1,\dots,\ivc t}\subseteq \ET(P)$, and denote the subposet $P\setminus\set{0_P}$ of $P$ by $Q$. Then the map
\begin{equation*}
g\colon \ET(P)\to \ET(Q),\qquad \vec x=(x_q:q\in P)\mapsto (x_q: q\in Q)
\end{equation*}
is well-defined. Indeed, if $\vec x$ satisfies the inequalities of \eqref{eq:elIg} for all $q\in P$, then these inequalities are also satisfied for all $q\in P\setminus\set{0_P}=Q$, and therefore $g(\vec x)\in\ET(Q)$. 

For each $q\in Q$, the upper covers of $q$ are the
same in $Q$ as in $P$. Hence, if
$\vec x,\vec y,\vec u,\vec v\in\ET(P)$ satisfy equations
\eqref{eq:ETPmeet} and \eqref{eq:ETPjoin}
for all $q\in P$, then these equations also
hold for all $q\in Q$.
That is, if $\vec u=\vec x\wedge \vec y$ in $\ET(P)$, then $g(\vec x\wedge \vec y)=g(\vec u)=g(\vec x)\wedge g(\vec y)$ in $\ET(Q)$. We obtain similarly that $g(\vec x\vee \vec y)=g(\vec v)=g(\vec x)\vee g(\vec y)$. Thus $g\colon \ET(P)\to \ET(Q)$ is a lattice homomorphism.

Let $\Psi:=\{(\vec x,\vec y)\in K^2:
g(\vec x)=g(\vec y)\}$; it is the congruence
kernel of $\restrict g K$.
By the simplicity of $K$, $\Psi$ is either
the equality relation or the full relation $K^2$.
If $\Psi$ were the equality relation, then
$\restrict g K$ would embed $K$ into $\ET(Q)$,
contradicting \eqref{eq:Pminsize}.
Hence $\Psi$ is the full relation.
Thus, the $t$ different $P$-tuples
$\ivc 1,\dots,\ivc t\in K$ have identical
components indexed by $q$ for every $q\in Q$.
Since these $P$-tuples are different,
their components indexed by $0_P$,
namely $x_{0_P}^{(1)},\dots,x_{0_P}^{(t)}$,
take $t\geq 3$ different values.
This contradicts \eqref{eq:elIg},
which allows at most two different values
for these coordinates, and completes the
proof of Proposition~\ref{prop:simple}.
\end{proof}

\end{document}